\documentclass[a4paper]{amsart}  

\usepackage{amssymb} \usepackage{amscd} \usepackage{times}

\def\zbb{\mathbb{Z}}  
  
  \def\phi{\varphi}
 \def\p1{{\mathbb{P}^1_\zbb}}

\begin{document}

\title{Asymptotic Estimate for Perturbed Scalar Curvature Equation.}
\author{Samy Skander Bahoura}

\address{6, rue Ferdinand Flocon, 75018 Paris, France. }
              
\email{samybahoura@yahoo.fr, bahoura@ccr.jussieu.fr} 

\date{}

\maketitle

\begin{abstract}

We consider the equation $ \Delta u_{\epsilon} = V_{\epsilon} {u_{\epsilon}}^{(n+2)/(n-2)}+ \epsilon W_{\epsilon}{u_{\epsilon}}^{\alpha } $ with $ \alpha \in ]\dfrac{n}{n-2}, \dfrac{n+2}{n-2} [ $ and we give some minimal conditions on $ \nabla V $ and $ \nabla W $ to have an uniform estimate for their solutions when $ \epsilon \to 0 $.
 
\end{abstract}

\bigskip

\bigskip

\begin{center}  1. INTRODUCTION AND RESULTS.
\end{center}

\bigskip

We denote $ \Delta=-\sum_i \partial_{ii} $ the geometric Laplacian on $ {\mathbb R}^n, n\geq 3 $.

\bigskip

Let us consider on open set $ \Omega $ of $ {\mathbb R}^n, n\geq 3 $, the following equation:

$$ \Delta u_{\epsilon} =V_{\epsilon} {u_{\epsilon}}^{(n+2)/(n-2)}+ \epsilon W_{\epsilon} {u_{\epsilon}}^{\alpha} \qquad (E_{\epsilon}) $$

where $ V_{\epsilon} $ and $ W_{\epsilon} $ are two regular functions and $ \alpha \in ] \dfrac{n}{n-2}, \dfrac{n+2}{n-2}[ $.

\bigskip

We assume:

$$ 0 < a \leq V_{\epsilon}(x) \leq b, \,\, ||\nabla V_{\epsilon}||_{L^{\infty}} \leq A \qquad (C_1) $$

$$ 0 < c \leq W_{\epsilon}(x) \leq d, \,\, ||\nabla W_{\epsilon}||_{L^{\infty}} \leq B \qquad (C_2) $$

{\bf Problem:} Can we have an $ \sup \times \inf $ estimate with the minimal conditions $ ( C_1) $ and $ ( C_2) $ ?

\bigskip

Note that for $ W \equiv 0 $, the equation $ (E_{\epsilon}) $  is the wellknowen scalar curvature equation on open set of $ {\mathbb R}^n $, $ n \geq 3 $. In this case, there is many results about this equation, see for example [B] and [C-L 1].

\bigskip

When $ \Omega = {\mathbb S}_n $ YY. Li, give a flatness condition to have the boundedness of the energy and the existence of the simple blow-up points, see [L1] and [L2].

\bigskip

In [C-L 2], Chen and Lin gave a conterexample of solutions of the scalar curvature equation with unbounded energy. The conditions of Li are minimal in heigh dimension.

\bigskip

Note that, in  [C-L 1] and [C-L 3], there is some results concerning Harnack inequalities of type  $ \sup \times \inf $ with the "Li-flatness" conditions for the following equation:

$$ \Delta u= Vu^{(n+2)/(n-2)}+ g(u) $$

where $ g $ is a regular function ( at least $ C^1 $ ) such that $ g(t)/[ t^{(n+2)/(n-2)}] $ is deacrising and tends to 0 when $ t \to +\infty $. They extend Li result ([L1]) to any open set of the euclidian space.

\bigskip

We can find in [A], some existence results for the presribed scalar curvature equation.

\bigskip

In our work we have no assumption on the energy. We use the blow-up analysis and the moving-plane method, developped by Gidas-Ni-Nirenberg, see [ G-N-N]. This method was used by different authors to have a priori estimates, look for example, [B], [B-L-S] ( in dimension 2), [C-L 1], [C-L 3], [L 1] and [L 2].

\bigskip

We set $ \delta =[(n+2)-\alpha (n-2)]/2 $, $ \delta \in ]0,1[ $. We have:

\bigskip

{ \it {\bf Theorem 1}. For all $ a,b,c,d, A, B >0 $, for all $ \alpha \in ]\dfrac{n}{n-2}, \dfrac{n+2}{n-2}[ $ and all compact set $ K $ of $ \Omega $, there is a positive constant $ c=c(a,b,c,d, A, B, \alpha, K, \Omega,n) $ such that:

$$ \epsilon^{(n-2)/2(1-\delta)} (\sup_K u_{\epsilon})^{1/3} \times \inf_{\Omega} u_{\epsilon} \leq c  $$

for all $ u_{\epsilon} $ solution of $ (E_{\epsilon}) $ with $ V_{\epsilon} $ and $ W_{\epsilon} $ satisfying the conditions $ (C_1) $ and $ (C_2) $.}

\bigskip

Now, we suppose that $ V_{\epsilon} $ satisfies:

$$ 0 < a \leq V_{\epsilon} (x) \leq b \,\,\, {\rm and} \,\,\, ||\nabla V_{\epsilon}||_{L^{\infty}(\Omega)} \leq k \epsilon \qquad (C_3) $$

We have:

{\it {\bf Theorem 2}. For all $ a,b,c,d, k, B >0 $, for all $ \alpha \in ]\dfrac{n}{n-2}, \dfrac{n+2}{n-2}[ $ and all compact set $ K $ of $ \Omega $, there is a positive constant $ c=c(a,b,c,d, k, B, \alpha, K, \Omega,n) $ such that:

$$ \sup_K u_{\epsilon} \times \inf_{\Omega} u_{\epsilon} \leq c  $$

for all $ u_{\epsilon} $ solution of $ (E_{\epsilon}) $ with $ V_{\epsilon} $ and $ W_{\epsilon} $ satisfying the conditions $ (C_3) $ and $ (C_2) $.} 

\bigskip

Note that in [B], we have some results as the previous but for prescribed scalar curvature equation with subcritical exponent tending to the critical. Here, we have a $ \sup \times \inf $ inequality for the scalar curvature equation, with critical exponent, perturbed by a nonlinear term. We can see the influence of this non-linear term.

\bigskip

\bigskip

\begin{center}

2. PROOFS OF THE THEOREMS.

\end{center}

\bigskip

\bigskip

\underbar {\bf Proof of the theorem 1.}

\bigskip

Without loss of generality, we suppose $ \Omega = B_1 $ the unit ball of $ {\mathbb R}^n $. We want to prove an a priori estimate around 0. We can also suppose $ \epsilon  \to 0 $, the case $ \epsilon \not \to 0 $ is solved in [B].

\bigskip

Let $ (u_i) $ and $ (V_i) $ be a sequences of functions on $ \Omega $ such that:

$$ \Delta u_i = V_i {u_i}^{(n+2)/(n-2)}+\epsilon_i W_i u_i^{\alpha }, \,\, u_i>0,  $$

with $ 0 < a \leq V_i(x) \leq b $, $ 0 < a \leq W_i(x) \leq d $, $ ||V_i||_{L^{\infty}} \leq A $ and $ ||W_i||_{L^{\infty}} \leq B $.

\bigskip

We argue by contradiction and we suppose that the $ \sup \times \inf $ is not bounded.

\bigskip

We have:

\bigskip

$ \forall \,\, c,R >0 \,\, \exists \,\, u_{c,R} $ solution of $ (E_1) $ such that:

$$ \epsilon^{(n-2)/2(1-\delta)} R^{n-2} (\sup_{B(0,R)} u_{\epsilon, c, R})^{1/3} \times \inf_{\Omega } u_{\epsilon, c, R} \geq c, \qquad (H) $$

\underbar {\bf Proposition :}{\it (blow-up analysis)} 

\smallskip

There is a sequence of points $ (y_i)_i $, $ y_i \to 0 $ and two sequences of positive real numbers $ (l_i)_i, (L_i)_i $, $ l_i \to 0 $, $ L_i \to +\infty $, such that if we set $ v_i(y)=\dfrac{u_i(y+y_i)}{u_i(y_i)} $, we have:

$$ 0 < v_i(y) \leq  \beta_i \leq 2^{(n-2)/2}, \,\, \beta_i \to 1. $$

$$  v_i(y)  \to \left ( \dfrac{1}{1+{|y|^2}} \right )^{(n-2)/2}, \,\, {\rm uniformly \,\, on \,\, all \,\, compact \,\, set \,\, of } \,\, {\mathbb R}^n . $$

$$ l_i^{(n-2)/2} {\epsilon_i}^{(n-2)/2(1-\delta)} [u_i(y_i)]^{1/3} \times \inf_{B_1} u_i \to +\infty,$$

\underbar {\bf Proof of the proposition:}

\bigskip

We use the hypothesis $ (H) $, we take two sequences $ R_i>0, R_i \to 0 $ and $ c_i \to +\infty $, such that,

$$ {\epsilon_i}^{(n-2)/2(1-\delta)} {R_i}^{(n-2)} (\sup_{B(0,R_i)} u_i)^{1/3} \times \inf_{B_1} u_i \geq c_i \to +\infty, $$

Let $ x_i \in  { B(x_0,R_i)} $ be a point such that $ \sup_{B(0,R_i)} u_i=u_i(x_i) $ and $ s_i(x)=(R_i-|x-x_i|)^{(n-2)/2} u_i(x), x\in B(x_i, R_i) $. Then, $ x_i \to 0 $.

\bigskip

We have:

$$ \max_{B(x_i,R_i)} s_i(x)=s_i(y_i) \geq s_i(x_i)={R_i}^{(n-2)/2} u_i(x_i)\geq \sqrt {c_i}  \to + \infty. $$ 

We set:

$$ l_i=R_i-|y_i-x_i|,\,\, \bar u_i(y)= u_i(y_i+y),\,\,  v_i(z)=\dfrac{u_i [y_i+\left ( z/[u_i(y_i)]^{2/(n-2)} \right )] } {u_i(y_i)}. $$

Clearly we have, $ y_i \to x_0 $. We also obtain:

$$ L_i= \dfrac{l_i}{(c_i)^{1/2(n-2)}} [u_i(y_i)]^{2/(n-2)}=\dfrac{[s_i(y_i)]^{2/(n-2)}}{c_i^{1/2(n-2)}}\geq \dfrac{c_i^{1/(n-2)}}{c_i^{1/2(n-2)}}=c_i^{1/2(n-2)}\to +\infty. $$

\bigskip

If $ |z|\leq L_i $, then $ y=[y_i+z/ [u_i(y_i)]^{2/(n-2)}] \in B(y_i,\delta_i l_i) $ with $ \delta_i=\dfrac{1}{(c_i)^{1/2(n-2)}} $ and $ |y-y_i| < R_i-|y_i-x_i| $, thus, $ |y- x_i| < R_i $ and, $ s_i(y)\leq s_i(y_i) $. We can write:

$$ u_i(y) (R_i-|y-y_i|)^{(n-2)/2} \leq u_i(y_i) (l_i)^{(n-2)/2}. $$

But, $ |y -y_i| \leq \delta_i l_i $, $ R_i >l_i$ and $ R_i-|y- y_i| \geq R_i-\delta_i l_i>l_i-\delta_i l_i=l_i(1-\delta_i) $. We obtain,

$$ 0 < v_i(z)=\dfrac{u_i(y)}{u_i(y_i)} \leq \left [ \dfrac{l_i}{l_i(1-\delta_i)} \right ]^{(n-2)/2}\leq 2^{(n-2)/2} . $$

We set, $ \beta_i=\left ( \dfrac{1}{1-\delta_i} \right )^{(n-2)/2} $, clearly, we have, $ \beta_i \to 1 $.

\bigskip

The function $ v_i $ satisfies:

$$ \Delta v_i= \tilde V_i {v_i}^{(n+2)/(n-2)}+\epsilon_i \tilde W_i \dfrac{v_i^{n/(n-2)}}{[u_i(y_i)]^{[(n+2)/(n-2)]-\alpha}} $$

where, $ \tilde V_i(y)=V_i \left [ y+y/[u_i(y_i)]^{2/(n-2)} \right ] $ and $ \tilde W_i(y)=W_i \left [ y+y/[u_i(y_i)]^{2/(n-2)} \right ] $.

\bigskip
 
Without loss of generality, we can suppose that $ \tilde V_i \to V(0)=n(n-2) $.

\bigskip

We use the elliptic estimates, Ascoli and Ladyzenskaya theorems to have the uniform convergence of $ (v_i) $  to $ v $ on compact set of $ {\mathbb R}^n $. The function $ v $ satisfies: 

$$ \Delta v=n(n-2)v^{N-1}, \,\, v(0)=1,\,\, 0 \leq v\leq 1\leq 2^{(n-2)/2}, $$

By the maximum principle, we have $ v>0 $ on $ {\mathbb R}^n $. If we use Caffarelli-Gidas-Spruck result, ( see [C-G-S]), we obtain, $ v(y)=\left ( \dfrac{1}{1+{|y|^2}} \right )^{(n-2)/2} $. We have the same properties that in [B].

\bigskip

\bigskip

\underbar {\bf Polar Coordinates} {\it (Moving-Plane method)}

\bigskip

Now, we must use the same method than in the Theorem 1 of [B]. We will use the moving-plane method.

\bigskip

We must prove the lemma 2 of [B].

\bigskip

We set $ t\in ]-\infty, -\log2 ] $ and $ \theta \in {\mathbb S}_{n-1}
$ :

$$ w_i(t,\theta)=e^{(n-2)t/2}u_i(y_i+e^t\theta), \,\,\,
\bar V_i(t,\theta)=V_i(y_i+e^t\theta)\,\, {\rm and} \,\, \bar W_i(t,\theta)=W_i(y_i+e^t \theta). $$

We consider the following operator $ L = \partial_{tt}-\Delta_{\sigma}-\dfrac{(n-2)^2}{4} $, with $
\Delta_{\sigma} $ the Laplace-Baltrami operator on $ {\mathbb
  S}_{n-1} $.

\bigskip

The function $ w_i $  is solution of:

$$ -Lw_i=\bar V_i{w_i}^{N-1}+\epsilon_i e^{[(n+2)-(n-2)\alpha] t/2} \bar W_i
{w_i}^{\alpha}. $$

For $ \lambda \leq 0 $ we set :

\bigskip

$ t^{\lambda}=2\lambda-t  $ $ w_i^{\lambda}(t,\theta)=w_i(t^{\lambda},\theta) $, $ \bar
V_i^{\lambda}(t,\theta)=\bar V_i(t^{\lambda},\theta) $ et $ \bar
W_i^{\lambda}(t,\theta)=\bar W_i(t^{\lambda},\theta) .$

\bigskip

\underbar {\bf Remark:} Here we work on $ [\lambda, t_i] \times {\mathbb S}_{n-1} $, with $ \lambda \leq -\dfrac{2}{n-2} \log u_i(y_i)+2 $ and  $ t_i \leq  \log \sqrt {l_i} $, where $ l_i $ is chooses as in the proposition.

\bigskip

First, like in [B], we have the following lemma:

\bigskip

\underbar{\bf Lemma 1:} 

\bigskip

Let $ A_{\lambda} $ be the following property:

\bigskip

$ A_{\lambda}=\{\lambda\leq 0,\,\,\exists \,\,
 (t_{\lambda},\theta_{\lambda })\in
   ]\lambda,t_i]\times  {\mathbb S}_{n-1},\,\,   {\bar
   w_i}^{\lambda}(t_{\lambda},\theta_{\lambda})-{\bar
   w_i}(t_{\lambda},\theta_{\lambda} ) \geq 0\} . $

\bigskip

Then, there is $ \nu \leq 0 $, such that for $\lambda \leq \nu $, $ A_{\lambda}  $ is not true.

\bigskip

Like in the proof of the Theorem 1 of [B], we want to prove the following lemma:

\underbar{\bf Lemma 2:} 

\bigskip

For $\lambda \leq 0 $ we have :

$$ {w_i}^{\lambda} -w_i \leq 0 \Rightarrow -L({w_i}^{\lambda} - w_i) \leq 0, $$

on $ ]\lambda,t_i]\times {\mathbb S}_{n-1} $.

\bigskip

Like in [B], we have:

\bigskip
  
\underbar{\bf A useful point: }

\bigskip

$ {\xi}_i= \sup \{\lambda \leq { \bar\lambda_i}=2+\log
  \eta_i, {w_i}^{\lambda} - w_i < 0 $, on $
  ]\lambda,t_i]\times {\mathbb S}_{n-1}
  \} $. The real $ \xi_i $ exists.

\bigskip

First, we have:

$$ w_i(2\xi_i-t,\theta)=w_i[(\xi_i-t+\xi_i-\log\eta_i-2)+(\log
\eta_i+2)] , $$

the definition of $ w_i $  and the fact that, $ \xi_i \leq t $, we obtain:

 $$ w_i(2\xi_i-t,\theta)=e^{[(n-2)(\xi_i-t+\xi_i-\log\eta_i-2)]/2}e^{n-2}v_i[\theta e^2e^{(\xi_i-t)+(\xi_i-\log\eta_i-2)}]
\leq 2^{(n-2)/2}e^{n-2}=\bar c. $$

\underbar {\bf Proof of the Lemma 2:}

\bigskip

We know that:

$$ -L( w_i^{\xi_i}-w_i)=[\bar V_i^{\xi_i
  }(w_i^{\xi_i})^{N-1}-\bar V_i
{w_i}^{N-1}]+ \epsilon_i [e^{\delta t^{\xi_i}}{\bar W_i}^{\xi_i
  }(w_i^{\xi_i})^{\alpha }-e^{\delta t}\bar W_i
{w_i}^{\alpha} ] ,  $$

with $ \delta =[(n+2)-(n-2)\alpha]/2 $.

\bigskip

We denote by $ Z_1 $ and $ Z_2 $ the following terms:

$$ Z_1=(\bar V_i^{\xi_i }-\bar V_i)(w_i^{\xi_i })^{N-1}+\bar
V_i[(w_i^{\xi_i })^{N-1}-{w_i}^{N-1}],$$

and

$$ Z_2=\epsilon_i (\bar
W_i^{\xi_i }-\bar W_i)(w_i^{\xi_i })^{\alpha }e^{\delta
  t^{\xi_i }}+ \epsilon_i e^{\delta
  t^{\xi_i }}\bar W_i[(w_i^{\xi_i })^{\alpha }-{w_i}^{\alpha
  }]+ \epsilon_i \bar W_i {w_i}^{\alpha }(e^{\delta t^{\xi_i }}-e^{\delta t} ). $$

But, using the same method as in the proof of the theorem 1 of [B], we have:

$$ {w_i}^{\xi_i} \leq w_i \,\,\,{\rm et } \,\,\,
w_i^{\xi_i}(t,\theta)\leq \bar c \,\,\, {\rm pour \, tout } \,\,\,
(t,\theta)\in [\xi_i,\log 2] \times {\mathbb S}_{n-1}  , $$

where  $ \bar c $ is a positive constantnot depending on  $ i
$ for $ \xi_i \leq \log \eta_i+2 $;

$$ |\bar V_i^{\xi_i }-\bar V_i|\leq A (e^t-e^{ t^{\xi_i }}) \,\,\,
{\rm et } \,\,\, |\bar W_i^{\xi_i }-\bar W_i|\leq B (e^t-e^{ t^{\xi_i
    }}), $$

Then,

\bigskip

$ Z_1 \leq A\, ({w_i^{\xi_i}})^{N-1} \,  (e^t-e^{ t^{\xi_i }}) \,\,\,
{\rm et} \,\,\,  Z_2 \leq  \epsilon_i B \, ({(w_i ^{\xi_i})}^{\alpha }\,  (e^t-e^{
  t^{\xi_i }})+ \epsilon_i c\, {(w_i^{\xi_i})}^{\alpha} \times   (e^{\delta t^{\xi_i
    }}-e^{\delta t} ) 
  $.

\bigskip

and,

\bigskip

$ -L(w_i^{\xi_i}-w_i) \leq (w_i^{\xi_i})^{\alpha}[ (A\,  {w_i^{\xi_i}}^{N-1-\alpha}+ \epsilon_i B)  \,  (e^t-e^{ t^{\xi_i }})+ \epsilon_i c\,\,  (e^{\delta t^{\xi_i
    }}-e^{\delta t} ) ].
$

\bigskip

But, $ w_i^{\xi_i} \leq \bar c $, we obtain:

\bigskip

$ -L(w_i^{\xi_i}-w_i)\leq  (w_i^{\xi_i})^{\alpha } [(A {\bar c}
^{N-1-\alpha}+ \epsilon_i B)  \,  (e^t-e^{ t^{\xi_i }})+ \epsilon_i c\,\,  (e^{\delta t^{\xi_i
    }}-e^{\delta t} ) ]. \,\,\,(1)$

\bigskip

We must see the sign of:

\bigskip

 $ \bar Z=[( A{\bar c}
^{N-1-\alpha}+ \epsilon_i B)  \, (e^t-e^{ t^{\xi_i }})+ \epsilon_i c\,\,  (e^{\delta t^{\xi_i
    }}-e^{\delta t} ) ]. $

\bigskip

But $ \alpha \in ] \dfrac{n}{n-2}, \dfrac{n+2}{n-2}[ $, $
\delta=\dfrac{n+2-(n-2)\alpha}{2} \in ]0,1[ $. 

\bigskip

For $ t \leq  t_i <0 $, we have:

$$ e^t \leq   e^{(1-\delta
  )t_i} e ^{\delta t} ,\,\,\, {\rm for \,\, all } \,\,\, t\leq t_i .$$

and the fact that $ t^{\xi_i}\leq t $ $(\xi_i \leq t )$, by integration of the previous two members, we obtain:

$$ e^t-e^{ t^{\xi_i }} \leq \dfrac{ e^{(1-\delta
  )t_i}}{\delta }(e^{\delta t}-e^{\delta t^{\xi_i
    }}), \,\,\, {\rm for\,\, all } \,\,\, t\leq t_i ,$$

We can write:

$$ (e^{\delta t^{\xi_i
    }}-e^{\delta t} ) \leq \dfrac{ \delta} { e^{(1-\delta
  )t_i} }\,  ( e^{ t^{\xi_i }}-e^t). $$

Then,

$$ -L(w_i^{\xi_i}-w_i) \leq (w_i^{\xi_i})^{\alpha}[-\dfrac{ \epsilon_i \delta \, c }{ e^{(1-\delta
  )t_i} }+ A \,{\bar c}^{N-1-\alpha}+\epsilon_i B]( e^t-e^{ t^{\xi_i }}). $$

The term $ \dfrac{\epsilon_i \delta \, c} {
  e^{(1-\delta
  )t_i}}  - A \,{\bar c}^{N-1-\alpha}-\epsilon_i B  $ is positive if:
   
   $$ \epsilon_i e^{-(1-\delta) t_i} \to +\infty, $$

then,

$$ \epsilon_i^{(n-2)/2(1-\delta)} e^{-(n-2)/2t_i} \to + \infty . $$ 
   
If we take, $ t_i= -\dfrac{2}{3(n-2)} \log u_i(y_i) $, we have:

$$ \epsilon_i^{(n-2)/2(1-\delta)} [u_i(y_i)]^{1/3}  \to + \infty. $$

It is given by our Hypothesis in the proposition.

\bigskip

But the Hopf Maximum principle, gives:

$$ \min_{\theta \in {\mathbb S}_{n-1}} w_i (t_i,\theta) \leq \max_{ \theta \in {\mathbb S}_{n-1}} w_i (2\xi_i - t_i,\theta), $$

then,

$$  e^{(n-2)t_i} u_i(y_i) \min_{B_2(0)} u_i \leq c, $$
   
and,

$$ [u_i(y_i)]^{1/3} \min_{B_2(0)} u_i \leq c, $$

Contradiction.   
   
\bigskip

\underbar {\bf Proof of the Theorem 2.}

\bigskip

The proof is similar than the proof of the theorem 1. Only the end of the proof is different.

\bigskip

\underbar {\bf Step 1:} The blow-up analysis give:

\bigskip

There is a sequence of points $ (y_i)_i $, $ y_i \to 0 $ and two sequences of positive real numbers $ (l_i)_i, (L_i)_i $, $ l_i \to 0 $, $ L_i \to +\infty $, such that if we set $ v_i(y)=\dfrac{u_i(y+y_i)}{u_i(y_i)} $, we have:

$$ 0 < v_i(y) \leq  \beta_i \leq 2^{(n-2)/2}, \,\, \beta_i \to 1. $$

$$  v_i(y)  \to \left ( \dfrac{1}{1+{|y|^2}} \right )^{(n-2)/2}, \,\, {\rm uniformly \,\, on \,\, all \,\, compact \,\, set \,\, of } \,\, {\mathbb R}^n . $$

$$ l_i^{(n-2)/2} u_i(y_i) \times \inf_{B_1} u_i \to +\infty,$$

\underbar {\bf Step 2:} Application of the Hopf maximum principle.

\bigskip

We have the same notation that in the proof of the theorem 1. First, we take $ t_i= \sqrt {l_i} $ as in the Step 1 and we look to the end of the proof of the theorem 1. We replace $ A $ by $ k \epsilon_i $. We want to proof that:

$$ {w_i}^{\lambda} -w_i \leq 0 \Rightarrow -L({w_i}^{\lambda} - w_i) \leq 0, $$

on $ ]\xi_i,t_i]\times {\mathbb S}_{n-1} $. We have the same defintion for $ \xi_i $ ( as in the proof of the theorem 1).

\bigskip

For $ t \leq  t_i <0 $, we have:

$$ e^t \leq   e^{(1-\delta
  )t_i} e ^{\delta t} ,\,\,\, {\rm for \,\, all } \,\,\, t\leq t_i .$$

and the fact that $ t^{\xi_i}\leq t $ $(\xi_i \leq t )$, by integration of the previous two members, we obtain:

$$ e^t-e^{ t^{\xi_i }} \leq \dfrac{ e^{(1-\delta
  )t_i}}{\delta }(e^{\delta t}-e^{\delta t^{\xi_i
    }}), \,\,\, {\rm for\,\, all } \,\,\, t\leq t_i ,$$

We can write:

$$ (e^{\delta t^{\xi_i
    }}-e^{\delta t} ) \leq \dfrac{ \delta} { e^{(1-\delta
  )t_i} }\,  ( e^{ t^{\xi_i }}-e^t). $$

Then,

$$ -L(w_i^{\xi_i}-w_i) \leq (w_i^{\xi_i})^{\alpha}[-\dfrac{ \epsilon_i \delta \, c }{ e^{(1-\delta
  )t_i} }+ k\epsilon_i \,{\bar c}^{N-1-\alpha}+\epsilon_i B]( e^t-e^{ t^{\xi_i }}). $$

The term $ \dfrac{\delta \, c} {
  e^{(1-\delta
  )t_i}}  - k  \,{\bar c}^{N-1-\alpha}- B  $ is positive because $ t_i \to -\infty $ and $ \delta \in ]0,1[ $.

\bigskip
  
But the Hopf Maximum principle, gives:

$$ \min_{\theta \in {\mathbb S}_{n-1}} w_i (t_i,\theta) \leq \max_{ \theta \in {\mathbb S}_{n-1}} w_i (2\xi_i - t_i,\theta), $$

then,

$$  e^{(n-2)t_i} u_i(y_i) \min_{B_2(0)} u_i \leq c, $$
   
and,

$$  l_i^{(n-2)/2} u_i(y_i) \min_{B_2(0)} u_i \leq c, $$

Contradiction with the step 1.   

\bigskip

\bigskip

\bigskip

\underbar {\bf R\'ef\'erences:}

\bigskip

[A] T. Aubin. Nonlinear Problems in Riemannian Geometry. Springer-Verlag 1998.

\smallskip

[B] S.S Bahoura. Majorations du type $ \sup u \times \inf u \leq c $ pour l'\'equation de la courbure scalaire prescrite sur un ouvert de $ {\mathbb R}^n, n\geq 3 $. J.Math.Pures Appl.(9) 83 (2004), no.9, 1109-1150.

\smallskip

[B-L-S] H. Brezis, Yy. Li Y-Y, I. Shafrir. A sup+inf inequality for some
nonlinear elliptic equations involving exponential
nonlinearities. J.Funct.Anal.115 (1993) 344-358.

\smallskip

[C-G-S] Caffarelli L, Gidas B., Spruck J. Asymptotic symmetry and local
behavior of semilinear elliptic equations with critical Sobolev
growth. Commun. Pure Appl. Math. 37 (1984) 369-402.

\smallskip

[C-L 1]  Chen C-C, Lin C-S. Estimates of the conformal scalar curvature
equation via the method of moving planes. Comm. Pure
Appl. Math. L(1997) 0971-1017.

\smallskip

[C-L 2] Chen C-C. and Lin C-S. Blowing up with infinite energy of conformal
metrics on $ {\mathbb S}_n $. Comm. Partial Differ Equations. 24 (5,6)
(1999) 785-799.

\smallskip

[C-L 3] Chen C-C, Lin C-S. Prescribing scalar curvature on $ {\mathbb S}_n $. I. A priori estimates.  J. Differential Geom.  57  (2001),  no. 1, 67--171.

\smallskip

[G-N-N] B. Gidas, W. Ni, L. Nirenberg, Symmetry and Related Propreties via the Maximum Principle, Comm. Math. Phys., vol 68, 1979, pp. 209-243.

\smallskip

[L1] Y.Y Li. Prescribing Scalar Curvature on $ {\mathbb S}_n $ and related Problems. I. J. Differential Equations 120 (1995), no. 2, 319-410.

\smallskip

[L2] Y.Y Li. Prescribing Scalar Curvature on $ {\mathbb S}_n $ and related Problems. II. Comm. Pure. Appl. Math. 49(1996), no.6, 541-597.

\end{document}